\documentclass[12pt]{amsart}
\usepackage{amsmath,amsthm,amsfonts,amssymb}
\begin{document} 
\newcommand{\supp}{\mathop{supp}}
\newcommand{\F}{{\mathbb F}}
\newcommand{\B}{{\mathbb B}}
\newcommand{\C}{{\mathbb C}}
\newcommand{\N}{{\mathbb N}}
\newcommand{\Q}{{\mathbb Q}}
\newcommand{\Z}{{\mathbb Z}}
\renewcommand{\P}{{\mathbb P}}
\newcommand{\R}{{\mathbb R}}
\newcommand{\rc}{\subset}
\newcommand{\tensor}{\otimes}
\newcommand{\rank}{\mathop{rank}}
\newcommand{\trace}{\mathop{tr}}
\newcommand{\dimc}{\mathop{dim}_{\C}}
\newcommand{\Lie}{\mathop{Lie}}
\newcommand{\Sing}{\mathop{Sing}}
\newcommand{\End}{\mathop{End}}
\newcommand{\Hol}{\mathop{Hol}}
\newcommand{\Auto}{\mathop{{\rm Aut}_{\mathcal O}}}
\newcommand{\alg}[1]{{\mathbf #1}}
\newtheorem{definition}{Definition}
\newtheorem*{claim}{Claim}
\newtheorem{corollary}{Corollary}
\newtheorem*{Conjecture}{Conjecture}
\newtheorem*{SpecAss}{Special Assumptions}
\newtheorem{example}{Example}
\newtheorem*{remark}{Remark}
\newtheorem*{observation}{Observation}
\newtheorem*{statement}{Statement}
\newtheorem*{fact}{Fact}
\newtheorem*{remarks}{Remarks}
\newtheorem{question}{Question}
\newtheorem{lemma}{Lemma}
\newtheorem{proposition}{Proposition}
\newtheorem{theorem}{Theorem}
\title{%
Non-degenerate Maps and Sets
}
%\shorttitle{}
\author {J\"org Winkelmann}
\begin{abstract}
We construct certain non-degenerate maps and sets,
mainly in the complex-analytic category.
\end{abstract}
%\subjclass{AMS Subject Classification:}
%
\address{%
J\"org Winkelmann \\
 Institut Elie Cartan (Math\'ematiques)\\
 Universit\'e Henri Poincar\'e Nancy 1\\
 B.P. 239, \\
 F-54506 Vand\oe uvre-les-Nancy Cedex,\\
 France
}
\email{jwinkel@member.ams.org\newline\indent{\itshape Webpage: }%
http://www.math.unibas.ch/\~{ }winkel/
}
\thanks{
{\em Acknowledgement.}
The author wants to thank
the Korea Institute for Advanced Study (KIAS) in Seoul.
The research for this article was partially done during the stay of the
author at this institute.}
\maketitle
\section{Introduction}
Let $X$ be a variety or an analytic space. A subset of resp.~a map to $X$
 is {\em degenerate}
if it resp.~its image is contained in a proper subvariety/analytic
subspace of $X$. Clearly it is easy to find degenerate maps and
subsets, but sometimes it is maybe not so easy to find
non-degenerate objects. In this paper our goal is to 
construct maps and infinite
subsets which are as non-degenerate as possible.

\section{Summary}
The main results of this article are the following:

Concerning the existence of non-degenerate maps and sets
in the complex-analytic category:
\begin{itemize}
\item
Let $X$ be an irreducible complex space and $S\subset X$ a countable subset.
Then there exists a holomorphic map $f$ from the unit disk $\Delta$
to $X$ with $S\subset f(\Delta)$.
\item 
For every irreducible complex space $X$ there exists a holomorphic
map from the unit disk $\Delta$ to $X$ with dense image.
More general: If $X$ and $Y$ are irreducible complex spaces
and $X$ admits a non-constant bounded holomorphic function,
then there exists a holomorphic map from $X$ to $Y$ with
dense image.
\item
Let $X$ be an irreducible complex space, $\dim(X)>0$. 
Then there exists an infinite subset
$\Gamma\subset X$ such that $\Gamma\cap Z$ is finite for every
closed analytic subset $Z\subsetneq X$.
\item
Every holomorphic map $f:\C\to\C^n$ can be approximated
uniformly on compact sets by a sequence $f_n$ of holomorphic
maps from $\C$ to $\C^n$ with dense image.
\end{itemize}
The corresponding statements in the real-analytic category are true
as well.

Furthermore, in the real-analytic category we have the following
statement (for which there is no complex-analytic analogue):

\begin{itemize}
\item
Let $X$ be a real analytic manifold, $\dim(X)>1$,
 and $D$, $D'$ discrete subsets.
Then every bijection $f_0:D\to D'$ extends to a real-analytic
diffeomorphism of $X$.
\end{itemize}

Finally we have an algebraic analogue of one of the complex-analytic
statements:
\begin{itemize}
\item
Let $X$ be a variety defined over a 
field $k$ such that the set of $k$-rational
points $X(k)$ is dense in the Zariski topology. Assume $\dim(X)>0$.

Then there exists an infinite subset $\Gamma\subset X(k)$
such that $\Gamma\cap Z$ is finite for every subvariety $Z\subsetneq X$.
\end{itemize}

\section{Preparations}
All manifolds are assumed to be paracompact, Hausdorff, finite-dimensional
and connected.

We start by recalling some standard facts on real analytic maps.

A crucial instrument for our purposes
 is Grauert's embedding theorem (see \cite{Gr}):
\begin{theorem}
Every real-analytic manifold can be embedded into some $\R^n$.
\end{theorem}
Here an \begin{em}embedding\end{em} is a proper injective immersion.

Using the standard embedding of $\R^n$ into the Stein complex manifold $\C^n$,
one obtains the following immediate consequence:
\begin{corollary}\label{emb-Stein}
Every real analytic manifold 
admits a real-analytic embedding 
 (as a closed real-analytic submanifold) 
into a Stein complex
manifold, namely into $\C^n$.
\end{corollary}
Grauert's embedding theorem is related to the following
approximation theorem (\cite{Gr}, see also \cite{H}):
\begin{theorem}
Let $X$, $Y$ be real-analytic manifolds.
Then $C^\omega_S(X,Y)$ is dense in $C^\infty_S(X,Y)$.
\end{theorem}

Here $C^\omega_S(X,Y)$ and $C^\infty_S(X,Y)$
are endowed with the Whitney topology (sometimes called {\em strong}
topology). 

Let us assume that $X$ and $Y$ are equipped with some Riemannian
metrics. 

Now a subbasis for this topology is given by the collection
of all sets 
\[
V(k,f,\epsilon)=\{F:X\to Y \text{ with }||D^k(F-f)(x)||<\epsilon(x)\ 
\forall x\in X\}
\]
where $k\in\N$, $f\in C^\infty(X,Y)$ and $\epsilon\in C^0(X,\R^+)$.
Here $D^k$ denotes the $k$-th derivative and the norm on $D^k$
is induced by the Riemannian metric.
See \cite{H} for more information on this topology.

In this article we will always use this topology if we deal with
smooth or real-analytic functions or maps. In contrast, we use
the topological of locally uniform convergence if we deal
with holomorphic maps.

The above approximation result admits a variant with interpolation
(\cite{T}, Thm~3.3, p.128):
\begin{theorem}\label{approx-inter}
Let $X$, $Y$ be real-analytic manifolds and $Z$ a closed analytic
subset of $X$.
Let $g:Z\to Y$ be a real-analytic map and $F$ denote the set of
all maps $f:X\to Y$ whose restriction to $Z$ coincides with $g$.

Then $C^\omega_S(X,Y)\cap F$ is dense in $C^\infty_S(X,Y)\cap F$.
\end{theorem}

We will use this result in the case where $Z$ is a discrete subset.

\begin{corollary}\label{ra-dense}
Let $M$ be a real analytic manifold.

Then there exists a real-analytic map $f:\R\to M$ with dense image.
\end{corollary}
\begin{proof}
Let $S\subset M$ be a countable dense subset of $M$ and
$\zeta:\N\to S$ a bijection. Now $\N$ is a discrete (hence
real-analytic)
subset of $\R$. Thus $\zeta:\N\to S\subset M$ extends to a
real-analytic
map $f:\R\to M$.
Finally $f(\R)$ is dense in $M$ because it contains $S$.
\end{proof}
\begin{lemma}\label{diffeo-open}
Let $X$ be a complete Riemannian manifold. Then the set of diffeomorphisms
is open in the set of self-maps $C^\infty(X,X)$
with respect to the Whitney topology.
\end{lemma}
\begin{proof}
See \cite{H}, Ch.2, Thm.1.7.
\end{proof}

\section{Automorphisms and discrete subsets}

\begin{theorem}\label{ext-bij}
Let $X$ be a real analytic manifold
of dimension at least two, $D$ and $D'$ discrete subsets.

Then every bijection $f_0:D\to D'$ extends to a real analytic
diffeomorphism $f$ of $X$.
\end{theorem}
\begin{proof}
Since $\dim(X)\ge 2$, the complement of a real curve is connected.
Therefore it is possible to find a family of disjoint smooth
curves $\rho_\gamma:[0,1]\to X$ ($\gamma\in D$) with $\rho_\gamma(0)=\gamma$
and $\rho_\gamma(1)=f_0(\gamma)$.
We choose disjoint open and relatively compact
neighbourhoods $U_\gamma$ of the curves 
$\rho_\gamma([0,1])$. Next on each $U_\gamma$ we choose a smooth
vector field $v_\gamma$ with support contained inside $U_\gamma$ such that
$(v_\gamma)_{\rho_\gamma(t)}=\gamma'_\gamma(t)$ for all $\gamma\in D$ and $t\in[0,1]$.
Then we define a global vector field $v$ by stipulating that
$v_x=(v_\gamma)_x$ if $x\in U_\gamma$ and $v_x=0$ if
$x\not\in\cup_\gamma 
U_\gamma$.
Note that the support $\supp(v)$ is contained in a disjoint
union of relatively compact open subsets. This implies that $v$
is globally integrable.
Now $\exp_v(1)$ is a smooth diffeomorphism of $X$ which extends $f_0$.
By theorem~\ref{approx-inter} there is a real-analytic map $f$
arbitrarily close to $\exp_v(1)$ such that $f|_D=f_0$.
 
Finally, since $\exp_v(1)$ is a
diffeomorphism and diffeomorphism are open in the Whitney topology,
we may require that $f$ is a diffeomorphism.
\end{proof}
\begin{remark}
If $\dim(X)=1$ (i.e.~$X\simeq\R$ or $X\simeq S^1$),
the statement of the theorem still holds, provided we
{\em assume} that $f_0$ extends to a homeomorphism of $X$.
\end{remark}

\begin{remark}
The corresponding statement for complex manifolds is wrong,
even for $X\simeq\C^n$.
Rosay and Rudin proved that there are infinite discrete subsets $D,D'$ in
$\C^n$ for all $n\in\N$ such that there exists no holomorphic
automorphism of $\C^n$ mapping $D$ to $D'$ (\cite{RR}).
In \cite{LDS} we generalized this result to the case
where $\C^n$ is replaced by an arbitrary Stein manifold.
\end{remark}
\begin{question}
Given a complex manifold $X$ ($\dim(X)>0$), is it always possible
to find discrete subsets $D,D'\subset X$ of the same cardinality
such that there exists no holomorphic automorphism $\phi$
of $X$ with $\phi(D)=D'$?
\end{question}
We believe that the answer is positive.
See \cite{LDS} for a more thorough discussion of this question.

\section{Maps: Disks with dense image}
\subsection{First Approach}
We show that every complex space admits a dense disc.

\begin{theorem}\label{disk}
Let $Z$ be an irreducible complex space
and let $S\subset Z$ be a countable subset (not necessarily discrete).

Then there exists a holomorphic map $F$ from the unit disk
$\Delta=\{z\in\C:|z|<1\}$ to $Z$ such that $S\subset F(\Delta)$.
\end{theorem}
\begin{proof}
Let $\sigma:\tilde Z\to Z$ be a desingularization 
(\cite{H},\cite{BM}) and $\tilde S\subset
\tilde Z$ a countable subset with $\sigma(\tilde S)=S$.
Without loss of generality we may assume that $S$ 
(hence also $\tilde S$) is an infinite subset.
Let $\zeta:\N\to\tilde S$ be a bijection.
Because $\tilde Z$ is connected, this map $\zeta$ extends
to a $C^\infty$-map $\zeta_1:\R\to\tilde Z$.
Using thm.~\ref{approx-inter}, it follows
that there is a real-analytic map $f_0:\R\to\tilde Z$
with $f_0|_{\N}=\zeta$.
Since analytic maps are locally given by convergent power series,
this map $f_0$ extends to a holomorphic map $f_1$ defined an some
open neighborhood $W$ of $\R$ in $\C$. 
Denoting the natural projection from $X=\C\times\tilde Z$
onto the second factor $\tilde Z$ by $pr_2$ we define
$F:W\to Z$ as $F=\sigma\circ pr_2\circ f_1$.
Note that $S\subset F(W)$.
By shrinking $W$, if necessary,
we may assume that $W\ne\C$ and that $W$ is simply-connected.
Then Riemann's mapping theorem (see e.g.~\cite{Rem})
implies that there is
a biholomorphic map $\Delta\simeq W$.
\end{proof}

\begin{corollary}\label{cor-dense-disk}
Let $X$ be an irreducible complex space.

Then there exists a holomorphic map $F$ from the unit disk
$\Delta=\{z\in\C:|z|<1\}$ to $X$ such that $\overline{F(\Delta)}=X$.
\end{corollary}
\begin{proof}
Choose a dense countable subset $S\subset X$ and invoke the preceding
theorem.
\end{proof}
\begin{remark}
In particular, it is always possible to connect two points in a complex
manifold by one disk.

Recall that for a complex manifold $X$ its 
{\em Kobayashi pseudodistance} $d_X$ is defined via chains of disks
(see \cite{K1},\cite{K2}). 
Since any two points can be connected by {\em one} disk,
one may also define a map $d':X\times X\to\R^+$ using just single
disks instead of disk chains. However, in this case it is not
clear whether $d'$ fulfills the triangle inequality. Thus it is
preferable to use disk chains.
\end{remark}
\begin{question}
Given a complex manifold $X$, is the set of all holomorphic maps
from the unit disk $\Delta$ to $X$ with dense image 
dense in the space of all holomorphic maps from $\Delta$ to $X$?
\end{question}
\begin{question}
Is there an analog of theorem~\ref{disk} for spaces over
non-archi\-me\-dean complete fields like e.g. the $p$-adic numbers?
\end{question}
\subsection{An alternative approach to dense disks}
By results of Fornaess and Stout (\cite{FS1},\cite{FS2})
it is known that for every complex
manifold $X$ of dimension $n$ there exists a surjective holomorphic map
from the polydisk $\Delta^n$ to $X$.
Thus in order to prove that there exists a holomorphic map from
the disk $\Delta$ to $X$ with dense image, it suffices to construct
a holomorphic map from $\Delta$ to $\Delta^n$ with dense image.
This can be done quite explicitly:
First we recall that $z\mapsto -i\frac{w+1}{w-1}$ 
maps $\Delta$ biholomorphically
to the upper half plane $H^+=\{z:\Im(z)>0\}$.

Now let $\lambda_1,\ldots,\lambda_{2n}$ be 
positive real numbers which are
linearly independent over $\Q$.
We consider the map
$
\phi:H^+\to\Delta^{2n}$ given
by 
\[
\phi:z\mapsto \left( e^{i\lambda_1z},\ldots,e^{i\lambda_{2n}z}\right)
\]
Fix a positive real number $\tau>0$.
Now for any $z\in H^+$ with $\Im(z)=\tau$ we have
$|e^{i\lambda_jz}|=e^{\Re(i\lambda_jz)}=e^{-\lambda_j\tau}$.
Thus the real line $i\tau+\R$ maps into the totally real $n$-dimensional
torus
given as
\[
\{(v_1,\ldots,v_n):|v_k|=e^{-\lambda_k\tau}\forall k\}.
\]
The linear independance of the $\lambda_j$ over $\Q$ ensures that
the image is dense.

Therefore the closure of $\phi(H^+)$ in $\Delta^{2n}$ consists of
all $v=(v_1,\ldots,v_j)$ such that either all $v_j$ are zero
or $\lambda_j\log|v_k|=\lambda_k\log|v_j|$ for all $j,k$.

As a next step we consider 
\[
\eta:v\mapsto \frac{1}{2}(v_1+v_2,v_3+v_4,v_5+v_6,\ldots).
\]
First, let us remark that given $r\ge s>0$ a complex number $z$
can be written as sum $z=z_1+z_2$ with $|z_1|=r$ and $|z_2|=s$
if and only if $r+s\ge|z|\ge r-s$.
Now observe that $\lim_{\tau\to\infty}e^{-\lambda\tau}=1$ for
$\lambda\in\R^+$.
As a consequence, for every
$(w_1,\ldots,w_n)\in\left(\Delta\setminus\{0\}\right)^n$
there is a number $\tau>>0$ such that
\[
\frac{1}{2}\left(e^{-\lambda_{2k}\tau}+e^{-\lambda_{2k-1}\tau}\right)
>|z_k|
>\frac{1}{2}\left(|e^{-\lambda_{2k}\tau}-e^{-\lambda_{2k-1}\tau}|\right)
\]
for all $k\in\{1,\ldots,n\}$
and consequently $\left(\Delta\setminus\{0\}\right)^n$ is contained
in the image $\eta(\overline{\Phi(H^+)})$.
Thus $\eta\circ\Phi:H^+\to\Delta^n$ has dense image.

In order to give a completely explicit example:
\[
z\mapsto\left(
\frac{1}{2}\left(e^{\frac{z+1}{z-1}}+e^{\sqrt{2}\frac{z+1}{z-1}}\right),
\frac{1}{2}\left(e^{\sqrt{3}\frac{z+1}{z-1}}+e^{\sqrt{5}\frac{z+1}{z-1}}
\right)
\right)
\]
defines holomorphic map from $\Delta$ to $\Delta^2$ with dense image.
\subsection{Other source manifolds}
Instead of considering maps from the unit disk we may investigate the
same question for other manifolds.
\begin{definition}\label{def-ud}
A connected complex manifold $X$ is called {\em universally dominating}
if for every irreducible complex space $Y$ there exists a holomorphic
map from $X$ to $Y$ with dense image.
\end{definition}

Taking $Y=\Delta$, it is clear that a universally dominating complex
manifold admits a non-constant bounded holomorphic function.
In fact this necessary condition is also sufficient.
\begin{theorem}
A connected complex manifold $X$ is universally dominating
(in the sense of def.~\ref{def-ud}) if and only if there exists
a non-constant bounded holomorphic function on $X$.
\end{theorem}
The crucial part of the proof is the lemma below.
\begin{lemma}
Let $U$ be a bounded connected open subset of $\C$.

Then there exists a holomorphic map $\rho:U\to\Delta$ with dense
image.
\end{lemma}
\begin{proof}
First we recall the theory of ``Ahlfors maps''.
(see \cite{Be}, p.~49.) %thm.13.1
If $\Omega$ is a connected bounded domain in $\C$
{\em with smooth boundary} 
(often called ``of finite type'') and $p\in\Omega$,
then there is a unique holomorphic map $\rho:\Omega\to\Delta$,
called ``Ahlfors map'' with the following properties:
\begin{enumerate}
\item
$\rho(p)=0$,
\item
$\rho'(p)\in\R^+$ and
\item $|\rho'(p)|\ge|f'(p)|$ for all holomorphic maps
$f:\Omega\to\Delta$ with $f(p)=0$.
\end{enumerate}
For a bounded domain $\Omega\subset\C$ with smooth boundary
this map $\rho$ is surjective and in fact even proper.

This theory implies in particular:

{\em Let $q\in\Delta$ with $\min\{1-|q|,|q|\}>\epsilon>0$ and
let $W=\{z\in\Delta:|z-q|>\epsilon\}$.

Then there exists a holomorphic map $\rho:W\to\Delta$ with
$|\rho'(0)|>1$.}

Indeed, if we take $p=0$ and compare $\rho$ with
$f=id_\Delta|_W$, then the unicity of the Ahlfors map $\rho$
implies
$|\rho'(p)|>|f'(0)|=1$.

Now let us consider an arbitrary bounded domain $U$, making no
assumptions about its boundary. Choose a basepoint $p\in U$.

Let ${\mathcal F}$ denote the family of all holomorphic maps
$f$ from $U$ to $\Delta$ with $f(p)=0$ and let
$\alpha=\sup_{f\in{\mathcal F}}|f'(p)|\in\R\cup\{+\infty\}$.
%\footnote{$*$}{Actually $\alpha<+\infty$ due to the lemma of
%Schwarz, but we do not use this fact.}
Let $f_n\in{\mathcal F}$ be a sequence with $\lim|f_n(p)|=\alpha$.
Using the theorem of Montel we may assume that the sequence of
holomorphic maps $f_n$ converges to a holomorphic map $\rho$.
By these arguments we see:

{\em For every bounded domain $U\subset\C$ and every point
$p\in U$ there exists a holomorphic map $\rho:U\to\Delta$ with
$\rho(p)=0$ such that $|\rho'(p)|\ge|f'(p)|$ for every holomorphic
map $f:U\to\Delta$ with $f(p)=0$.}

Unlike an Ahlfors map for a bounded domain with {\em smooth boundary},
in general this map $\rho$ is not necessarily surjective.
For instance, consider $U=\Delta\setminus\{1/2\}$ and $p=0$.
In this case we have $\rho=\lambda id_\Delta|_U$ for some 
$\lambda\in\C$ with $|\lambda|=1$.

However, we claim that $\rho$ has always dense image.
Indeed, let us assume that the image is not dense.
Then we can find a relatively compact ball $B_\epsilon(q)\subset \Delta$
such that $B_\epsilon(q)\cap \rho(U)=\emptyset$.
But then we obtain a contradiction to the maximality property
of $\rho$ by considering the composition of $\rho$ with
the Ahlfors map of $\Delta\setminus\overline{B_\epsilon(q)}$.
\end{proof}

\begin{proof}[Proof of the theorem]
Let $f$ be a non-constant bounded holomorphic function on $X$.
The image $f(X)$ is a bounded open subset of $\C$.
By the preceding lemma there is a holomorphic map $\rho$
from $f(X)$ to $\Delta$ with dense image.
Thus there is a holomorphic map from $X$ to $\Delta$ with dense
image, namely $\rho\circ f$.

Let $Y$ be an irreducible complex space. Thanks to cor.~\ref{cor-dense-disk}
there is a holomorphic map $F:\Delta\to Y$ with dense image.
Then $F\circ\rho\circ f:X\to Y$ is a holomorphic map with dense
image.
\end{proof}

\section{Entire Curves}
\subsection{Entire Curves in the Affine Space}
The purpose of this section is to prove the following statement.
\begin{proposition}\label{prop-ent-curve}
For a holomorphic map $f:\C\to\C^n$ let $J_d(f):\C\to\C^{nk}$
be defined as $J_d(f)=(f,f',f'',\ldots)$.

Let $d\in\N\cup\{0\}$ and $\Phi_d\subset\Hol(\C,\C^n)$
the set of all holomorphic maps from $\C$ to $\C^n$ for which
$\overline{(J_df)(\C)}=\C^{nk}$.

Then $\Phi_d$ is dense in $\Hol(\C,\C^n)$ with respect to the topology
of locally uniform convergence.
\end{proposition}

To prepare the proof of the proposition, we develop some lemmata.
\begin{lemma}
Let $p\in\C$ and $|p|>r>0$ and $\epsilon>0$.

Then there exists a polynomial $P\in\C[X]$ with $P(p)=1$
and $|P(z)|\le\epsilon$ for all $z\in B_r=\{z\in \C:|z|\le r\}$.
\end{lemma}
\begin{proof}
Choose $N\in\N$ such that $(\frac{r}{|p|})^N<\epsilon$
and take $P(z)=\left(\frac{z}{p}\right)^N$.
\end{proof}

\begin{lemma}\label{lemma3}
Let $p\in\C$, $|p|>r>0$, $\epsilon>0$, $d\in\N$ and
$(a_0,\ldots,a_n)\in\C^{n+1}$.

Then there exists a polynomial $P\in\C[X]$ such that
\begin{enumerate}
\item
$|P(z)|<\epsilon$ whenever $|z|\le r$ and
\item
$P^{(k)}(p)=a_k$ for $k=0,\ldots,d$.
\end{enumerate}
\end{lemma}
(Here $P^{(k)}$ denotes the $k$-th derivative of $P$.)
\begin{proof}
We define recursively polynomials $P_m$ ($m=0,\ldots,d$)
such that
\begin{enumerate}
\item
$|P_m(z)|<\frac{m+1}{d+1}\epsilon$ whenever $|z|\le r$ and
\item
$P^{(k)}(p)=a_k$ for $k=0,\ldots,m$.
\end{enumerate}
The existence of $P_0$ follows from the preceding lemma.
Now let us assume that $P_{m-1}$ has been constructed.
Define
\[
Q_m(z)=\left( a_m-P^{(m)}_{m-1}(p)\right)\frac{1}{m!}(z-p)^m
\]
and choose $R_m$ such that $R_m(p)=1$ and 
$|R_m(z)|<\frac{\epsilon}{d+1}\inf_{|w|<r}|Q_m(w)|^{-1}$ 
for all $z\in\Delta_r$.

Then $P_m=P_{m-1}+Q_mR_m$ has the desired properties.
Finally set $P=P_d$.
\end{proof}
\begin{lemma}\label{cs}
Let $D\subset\C$ be a discrete subset with $\min\{|z|:z\in D\}>r>0$,
$d\in\N\cup\{0\}$, $f:D\to \C^{d+1}$ a map and $\epsilon>0$.

Then there exists a holomorphic function $F$ on $\C$ with
$J_dF(\gamma)=f(\gamma)$ for all $\gamma\in D$ and
\[
\max_{|z|\le r}|F(z)|\le\epsilon
\]
\end{lemma}
\begin{proof}
Fix an enumeration $n\mapsto\gamma_n$ of $D$
such that $|\gamma_{n+1}|\ge|\gamma_n|$ for all $n\in\N$.
Choose a strictly increasing sequence $r_n$ with $r<r_n<|\gamma_n|$
and $\lim r_n=+\infty$.
(E.g.~choose $r_1=(r+|\gamma_1|)/2$ and
$r_{n+1}=(r_n+|\gamma_n|)/2$ for $n\ge 1$.)
Now we define recursively a sequence of functions $F_n$ as follows:
$F_0=0$. Assume $F_1,\ldots, F_{n-1}$ are already defined.
Let $Q_n$ be a polynomial with $(J_dQ_n)(\gamma_j)=0$ for $j<n$
and $(J_dQ_n)(\gamma_n)=f(\gamma_n)-(J_dF_{n-1})(\gamma_n)$.
Choose $\delta$ such that $\delta||Q_n||_{B_{r_n}}<2^{-n}\epsilon$.
Due to lemma~\ref{lemma3} there is a polynomial $P_n$ with
$P_n(\gamma_n)=1$,
$P_n^{(k)}(\gamma_n)=0$ for $1\le k \le d$ 
and $||P_n||_{B_{r_n}}<\delta$.
Now define $F_n(z)=F_{n-1}(z)+P_n(z)Q_n(z)$.
Observe that
\[
(J_dF_n)(\gamma_n)=(J_dF_{n-1})(\gamma_n)+(J_d(P_nQ_n))(\gamma_n)
%=F_{n-1}(\gamma_n)+1\left(f(\gamma_n)-F_{n-1}(\gamma_n)\right)
=f(\gamma_n)
\]
and
\[
(J_dF_n)(\gamma_j)=(J_dF_{n-1})(\gamma_j)+(J_d(P_nQ_n))(\gamma_j)
=(J_dF_{n-1})(\gamma_j)
\]
for $j<n$.
Hence $(J_dF_n)(\gamma_j)=f(\gamma_j)$ for all $j\le n$ by induction.

Since $||P_nQ_n||_{B_{r_n}}<2^{-n}\epsilon$, the sequence of functions $F_n$
converges to a global holomorphic function $F$. By construction $F$
has the desired properties.
\end{proof}

Now we can prove the proposition.
\begin{proof}[Proof of the proposition]
We have to show that for every $R,\epsilon>0$, $f\in\Hol(\C,\C^n)$
there exists a map $\phi\in\Phi$ with $||\phi-f||_{B_R}<\epsilon$.

Let $D=\{n\in\N:n>R\}$ and fix a bijection $\xi:D\to(\Q[i])^n$.
By lemma~\ref{cs} there is a holomorphic map $g:\C\to\C^n$
with $g(n)=\xi(n)-f(n)$ for all $n\in D$ and
$||g||_{B_R}<\epsilon$.
Now $\phi=g+f$ is a map from $\C$ to $\C^n$ with dense image such that
$||\phi-f||_{B_R}<\epsilon$.
\end{proof}

\subsection{Entire Curves in Projective Varieties}
Next we discuss dense entire curves in certain projective varieties.
As a preparation we prove the lemma below.

\begin{lemma}\label{lemma-Z}
Let $Z\subset\C^n$ be an algebraic subvariety of codimension at least
$2$. Then there exists a dominant morphism $F:\C^n\to\C^n\setminus Z$.
\end{lemma}
\begin{proof}
Let $e_1,\ldots,e_n$ denote the standard basis of $\C^n$,
$H_i$ the quotient vector space $\C^n/\C e_i$ and $\pi:\C^n\to H_i$
the natural linear projection. Let $Z_i$ denote the closure of
$\pi_i(Z)$.
Then $Z_i$ is an algebraic subvariety of codimension at least one.

Without loss of generality we may assume that none of the $Z_i$
contains the origin $0\in H_i$. Now we choose polynomial functions
$P_i$ on $H_i$ such that $P_i|_{Z_i}\equiv 0$ and $P_i(0)=1$.
For each $i\in\{1,\ldots,n\}$ and $\zeta\in\C$
we obtain an automorphism $\phi_{i,\zeta}:v\mapsto
v+\zeta P_i(\pi_i(v))e_i$. Now for $t=(t_1,\ldots,t_n)\in\C^n$
we define an automorphism $\Phi_t$ of $\C^n$ by
\[
\Phi_t:v\mapsto \phi_{n,t_n}\circ\ldots\circ\phi_{1,t_1}(v)
\]
By the construction each of the $\phi_{i,\zeta}$ stabilizes $Z$
pointwise.
Combined with $0\not\in Z$ and the invertibility of $\Phi_t$
we obtain that $\Phi_t(0)$ is never contained in $Z$.
Therefore we obtain an algebraic morphism $F:\C^n\to\C^n\setminus Z$
by
\[
F(t)=\Phi_t(0).
\]
The definition of $\Phi$ also implies that $F(\zeta e_i)=\zeta e_i$
for every $i\in\{1,\ldots,n\}$, $\zeta\in\C$.
It follows that $DF$ has maximal rank at $0$. Thus $F$ must be
dominant.
\end{proof}

\begin{proposition}\label{prop-unirational}
Let $X$ be a unirational complex projective variety.

Then there exists a holomorphic map $f:\C\to X$ with dense
image.
\end{proposition}
\begin{proof}
By the definition of unirationality we obtain a dominant rational
map $f_0:\C^n\to X$. The indeterminacy locus $E(f_0)$ is an algebraic
subvariety of codimension at least two and $f_0$ is regular outside
$E(f_0)$.
By lemma~\ref{lemma-Z} there is a dominant morphism
$g:\C^n\to\C^n\setminus
E(f_0)$. Thus we obtain a dominant morphism $h:\C^n\to X$ via
$h=f_0\circ g$.
By prop.~\ref{prop-ent-curve} there is a holomorphic map
$\phi:\C\to\C^n$
with dense image. Therefore $f=h\circ\phi:\C\to X$ is a holomorphic
map with dense image.
\end{proof}

\begin{question}
Let $X$ be a complex projective variety. Is $X$ {\em special} in the
sense of Campana (\cite{Cam}) iff there exists a holomorphic
map $f:\C\to X$ with dense image?
\end{question}

A positive answer would in particular require that there is such a
map $f:\C\to X$ for every $K3$-surface $X$. By the results of
Buzzard and Lu (\cite{BL}) there is a holomorphic map $F:\C^2\to X$ with dense
image if $X$ is an elliptic $K3$-surface or a Kummer surface.
Combined with prop.~\ref{prop-ent-curve}
 this implies that there is a holomorphic map
$f:\C\to X$ with dense image for these two special kinds of
$K3$-surfaces. However, for arbitrary $K3$-surfaces this is still
an open question.
\subsection{Other ground fields}
In this section we have rarely used special properties of the
field of complex numbers. Thus the results mostly are valid
for arbitrary fields with an absolute value
(like e.g. the field $\Q_p$ of $p$-adic numbers, or $\Omega_p$
or the field $\F_q((t))$ 
of formal Laurent series over some finite field $\F_q$.)

However, the following restrictions are necessary:
\begin{itemize}
\item
If the ground field is not locally compact, ``uniform convergence
on compact sets'' must be replaced by ``uniform convergence
on bounded sets'' and the condition ``$D$ is discrete in $\C$''
has to be replaced by ``{\em $\{z\in D:|z|<R\}$ is finite for all $R>0$}''.
\item
If the characteristic $p$ of the field is positive, one cannot divide
by $p$ and the $p$-th derivative is always zero.
Hence prop.~\ref{prop-ent-curve} and the lemmata \ref{lemma3}
and \ref{cs}
remain valid
only in the case where $d$ is smaller than the characteristic.
\item
A dominant morphism between algebraic varieties has necessarily
dense image
only if the ground field is algebraically closed.
Thus one needs to assume that the ground field is algebraically
closed for prop.~\ref{prop-unirational} and lemma~\ref{lemma-Z}.
\end{itemize}

\section{Non degenerate Sets}
\subsection{Analytic case}

\begin{proposition}
Let $X$ be an irreducible complex (resp.~real) analytic space.
Then there exists an infinite subset $\Gamma\subset X$ such that
$Z\cap\Gamma$ is finite for every closed analytic subset $Z$ of $X$.
\end{proposition}
\begin{proof}
It suffices to consider the real case.

We need some preparations.
First we recall that
\[
\arctan:\R\to ]-1,+1[
\]
is a proper real analytic map.

Next we choose a point $x\in X\setminus\Sing(X)$.
Define $D_r=\{x\in\R^d:||x||<r\}$.
Fix a real-analytic coordinate chart on an open neighbourhood
$U$ of $x$ in $X\setminus\Sing(X)$ as
\[
\zeta: U(x) \stackrel{\sim  }{\longrightarrow}D_2
\]
Let $d=\dim(X)$.
Due to cor.~\ref{ra-dense}
 there exists a real-analytic map 
\[
\phi:\R\to\partial D_1=\{x\in\R^d:||x||=1\}
\]
with dense image.
Let $F:\R\to U\subset X$ be the real-analytic map given by
\[
F:t\mapsto \zeta^{-1}\left( \arctan(t)\phi(t)\right)
\]
We claim that  $\overline{F(\R)}\supset\zeta^{-1}(\partial D_1)$.
Indeed, for every $x\in \partial D_1$ there is a
sequence $t_n\in\R$ with $\lim t_n=+\infty$ and $\lim\phi(t_n)=x$,
since $\phi(\R)$ is dense in $\partial D_1$. This implies $\lim F(t_n)=x$,
because $\lim\arctan(t_n)=1$. On the other hand $F(\R)\cap
\zeta^{-1}(\partial D_1)=\{\}$.

Let $\Gamma_0$ be a bounded infinite subset of $\R$
and define our infinite subset $\Gamma$ as
\[
\Gamma= F(\Gamma_0).
\]

In order to prove the proposition, we have to show that every
infinite subset of $\Gamma$ is Zariski dense in $X$.

Let $\Gamma'$ be an infinite subset of $\Gamma$ and $A$ the smallest
closed real analytic subset of $X$ containing $\Gamma'$.
We want to show that $A=X$.
Because $\Gamma$ is contained in the compact set $\zeta^{-1}(\overline{D_1})$
the analytic set $A$ has only finitely many irreducible components.
Hence there is no loss in generality in assuming that $A$ is
irreducible.
Since $\Gamma_0$ is bounded, there is an accumulation point $q$ of
$\Gamma'$ with $||\zeta(q)||<1$. 
Thus $F^{-1}(\Gamma')\subset F^{-1}(A)$
contains an accumulation point $p=F^{-1}(q)$.
By the identity
principle this implies that $F(\R)\subset A$.
Therefore $\overline{F(\R)}\subset A$. However, $\overline{F(\R)}$
contains
$C=\zeta^{-1}(\partial D_1)$ and $C$ is a real-analytic hypersurface.
Because of the irreducibility of $A$ we thus obtain that either
$A=C$ or $A=X$. But $A=C$ is impossible, since $C$ does not
contain any element of $\Gamma$. Hence $A=X$.
Thus every infinite subset $\Gamma'$ of $\Gamma$ is
dense in $X$ with respect to the analytic Zariski topology.
It follows that for every proper closed analytic subset $Z\subsetneq X$
the intersection $Z\cap\Gamma$ must be finite.
\end{proof}
\subsection{Algebraic case}
\begin{proposition}
Let $V$ be a variety defined over a field $k$.
Assume that $V(k)$ is Zariski-dense.

Then there exists an infinite subset $\Gamma\subset V(k)$ such that
for every proper subvariety $Z$ of $V$ the intersection $\Gamma\cap Z$
is finite.
\end{proposition}
\begin{proof}
Let $Div^+(V)$ the set of all effective reduced $k$-Weil divisors
on $V$. Using Chow schemes, $Div^+(V)$ can be exhausted
by an increasing sequence of $k$-varieties $C_n$. We define the universal space
$U_{n,m}\subset C_n\times V^m$ by the condition that $(\theta;
x_1,\ldots,x_m)\in U_{n,m}$  iff all the points $x_i$ are contained
in the support of the Weil divisor $D_\theta$ indexed by $\theta$.
Note that $U_{n,m}$ is of codimension at least $m$.
%(Here and in the sequel the dimension of the empty set is $-\infty$ and
%consequently the codimension of the empty set is $+\infty$.)
Let $\rho_{n,m}:U_{n,m}\to V^m$ be the natural projection.
Its generic fiber dimension does not exceed $\dim(C_n)-m$.
(Here the dimension of the empty space is defined as $-\infty$.)
We define $\Omega_{n,m}\subset V^m(k)$ as the set of 
$k$-rational points where
the fiber dimension of $\rho_{n,m}$ is less or equal $\dim(C_n)-m$.
Thus for $m>\dim(C_n)$ the fiber $\rho_{n,m}^{-1}(S)$ is empty
for every $S\in\Omega_{n,m}$.
Observe that $\Omega_{n,m}$ contains a Zariski open subset of $V^m(k)$.
Let $S=(s_1,\ldots,s_m)\in\Omega_{n,m}$.
We may choose a point in each irreducible component of
the fiber $\rho^{-1}_{n,m}(S)$. Let $\Theta_1,\ldots,\Theta_l$
denote the corresponding Weil divisors. Then for every 
$x\in V\setminus\cup_{i=1}^l\,\Theta_i$ we have
\[
(S,x)=(s_1,\ldots,s_m,x)\in\Omega_{n,m+1},
\]
because the fiber $\rho_{n,m+1}^{-1}(S,x)$ is a subvariety
of $\rho_{n,m}^{-1}(S)$ which by construction does not contain
any of the $(\Theta_i,S,x)$.
In particular, if $S\in\Omega_{n,m}$, then there exists a Zariski
open subset $W$ of $V$ such that $(S,x)\in\Omega_{n,m+1}$ for every
$x\in W(k)$.
Moreover $W(k)\ne\{\}$, because $V(k)$ is assumed to be Zariski dense
in $V$.

Therefore it is possible to choose recursively
a sequence $\gamma_j\in V(k)$
with the following property:

{\em If $n\le r\le j$ and if $S$ is a 
finite subset of $\{\gamma_r,\ldots,\gamma_j\}$,
then $S\in\Omega_{n,\#S}$.}

(If $\gamma_1,\ldots,\gamma_{j-1}$ are already chosen, this property
defines an open subset of $V$ from which we have to choose $\gamma_j$.)

Now this property is equivalent to the following statement:

{\em
Let $S$ be a finite subset of $\Gamma=\{\gamma_j:j\in\N\}$
with $r=\min\{j:\gamma_j\in S\}$.

Then $S\in\Omega_{n,\#S}$ for all $n\le r$.}

Recall that for $m>\dim(C_n)$ the fiber $\rho^{-1}_{n,m}(S)$ 
is empty for all $S\in\Omega_{n,m}$.
It follows that if $S$ is a finite subset of $\Gamma$
and $r=\min\{j:\gamma_j\in S\}$, then
$S$ is not contained in the support for any Weil divisor 
parametrized by $\theta\in C_n$ for any $n$ with $n\le r$ and
$\#(S)>\dim(C_n)$.

Let $Z$ be a proper subvariety of $V$. Then $Z$ is contained in
a Weil divisor parametrized by an element $\theta\in C_n$ for some
$n$.
Now the above arguments imply that the cardinality of a finite
subset $S$ of $(\Gamma\cap Z)\setminus\{\gamma_1,\ldots,\gamma_{n-1}\}$
is bounded by $\dim(C_n)$.
Hence $\#(\Gamma\cap Z)\le \dim(C_n)+n-1$.
\end{proof}
\begin{remark}
For $k=\C$ this statement can also be deduced from the analytic
analogue stated above. However, already for proper subfields
$k\subsetneq \C$ this is not just a corollary, because
in this case the analytic statement only yields the existence of
such a subset $\Gamma$ in $V(\C)$, but not that $\Gamma$ can be
chosen inside $V(k)$.
\end{remark}

\end{document}